\documentclass[11pt,reqno]{amsart}

\usepackage{amssymb,amsmath,epsfig}
\usepackage{amsfonts}
\usepackage{pst-plot}

\newcommand{\Z}{\mathbb{Z}}
\newcommand{\R}{\mathbb{R}}
\newcommand{\C}{\mathbb{C}}

\def\H{\mathbb{H}}
\newcommand{\leg}[2]{\genfrac{(}{)}{}{}{#1}{#2}}
\newtheorem{theorem}{Theorem}[section]

\newtheorem{lemma}[theorem]{Lemma}
\newtheorem{corollary}[theorem]{Corollary}

\newtheorem*{theorem*}{Theorem}
\newtheorem{remark}[theorem]{Remark}

\numberwithin{equation}{section}

\title[Generalized rank moments and Durfee symbols]
{Automorphic properties of generating functions for generalized rank moments and Durfee symbols}

\date{\today}
\author{Kathrin Bringmann, Jeremy Lovejoy, and Robert Osburn}
\address{School of Mathematics, University of Minnesota, Minneapolis, MN 55455, U. S. A.}

\address{CNRS, LIAFA, Universit\'e Denis Diderot,
2, Place Jussieu, Case 7014, F-75251 Paris Cedex 05, FRANCE}

\address{School of Mathematical Sciences, University College Dublin, Belfield, Dublin 4, Ireland}

\email{brigmann@math.umn.edu}

\email{lovejoy@liafa.jussieu.fr}

\email{robert.osburn@ucd.ie}

\thanks{The second author was partially supported by an ACI
``Jeunes Chercheurs et Jeunes Chercheuses".  The second and third
authors were partially supported by a PHC Ulysses}


\begin{document}
\begin{abstract}
We define two-parameter generalizations of two combinatorial
constructions of Andrews:  the $k$th symmetrized rank moment and
the $k$-marked Durfee symbol.  We prove that three specializations
of the associated generating functions are so-called quasimock theta
functions, while a fourth specialization gives quasimodular forms.
We then define a two-parameter generalization of Andrews' smallest
parts function and note that this leads to quasimock theta
functions as well.  The automorphic properties are deduced using
$q$-series identities relating the relevant generating functions
to known mock theta functions.
\end{abstract}

\maketitle

\section{Introduction}
The series $\mathcal{N}_{2v}(0,0;q)$, defined for $v \geq 1$ by
\begin{equation*}  
\mathcal{N}_{2v}(0,0;q) :=
\prod_{k = 1}^{\infty}\frac{1}{1-q^k}\sum_{n \in \mathbb{Z}
\setminus\{ 0\}} \frac{(-1)^{n-1}q^{n(3n-1)/2+vn}}{\left(1-q^n \right)^{2v}},
\end{equation*}
have recently been the focus of several studies.  These series are
tied to some new and exciting partition-theoretic constructions
\cite{An1} and their coefficients satisfy many elegant identities,
congruences, and asymptotic properties \cite{An1,Br1,Br-Ga-Ma1}.
From a theoretical standpoint, such results are ultimately due to
the fact that the series $q^{-1}\mathcal{N}_{2v}(0,0;q^{24})$
enjoy the automorphic structure of \emph{quasimock theta
functions} (see below for the definition).


Quasimock theta functions were introduced in \cite{Br1, Br-Ga-Ma1}
to describe functions like $\mathcal{N}_{2v}(0,0;q)$ which
resemble Ramanujan's mock theta functions but involve additional
quasimodular components.
In the same way that the modularity of the generating function for
partitions has many important consequences, the theory of
quasimock theta functions can be used to prove important
properties of the combinatorial objects encoded in their
coefficients.

This paper lies at the intersection of two questions. First, are
there meaningful generalizations of the combinatorial
constructions associated with $\mathcal{N}_{2v}(0,0;q)$?  Second,
how can one find further examples of simple $q$-series which are
quasimock theta functions?  The answer to both of these questions
will be found in the series $\mathcal{N}_{2v}(d,e;q)$, defined for
$v \geq 1$ by
\begin{equation}
\mathcal{N}_{2v}(d,e;q)
:= \frac{(-dq,-eq)_{\infty}}{(q,deq)_{\infty}}\sum_{n \in
\mathbb{Z} \setminus \{0\}} \frac{(-1)^{n-1}q^{{n+1 \choose 2} +
vn}(de)^n\left(-1/d,-1/e\right)_n}{(1-q^n)^{2v} \left(-dq,-eq\right)_n}.
\label{eta2vofrsn}
\end{equation}
Here we employ the standard $q$-series notation,
$$
\left(a_1,a_2,\dots,a_j;q \right)_n :=
\frac{\left(a_1,a_2,\dots,a_j;q\right)_{\infty}}{\left(a_1q^n,a_2q^n,\dots,a_jq^n;q\right)_{\infty}},
$$
where
$$
(a_1,a_2,\dots,a_j;q)_{\infty} :=
\prod_{k=0}^{\infty}\left(1-a_1q^k\right)\left(1-a_2q^k\right)\cdots\left(1-a_jq^k\right).
$$
We follow the custom of dropping the $;q$ unless the base is
something other than $q$.

In the first part of the paper we discuss two combinatorial
interpretations of the series $\mathcal{N}_{2v}(d,e;q)$. The
natural context is that of overpartition pairs
\cite{Br-Lo1,Lo1.5,Lo-Ma1}. In this setting we shall define a
\emph{generalized $k$th symmetrized rank moment} and a
\emph{generalized $k$-marked Durfee symbol}, each of which has
$\mathcal{N}_{2v}(d,e;q)$ as its generating function (see Section
2). When $d=e=0$ we recover the partition-theoretic work of
Andrews \cite{An1}.

In the second part of the paper we show that the series
$\mathcal{N}_{2v}(1,0;q)$, $\mathcal{N}_{2v}(1,1/q;q^2)$, and\\
$q^{-1}\mathcal{N}_{2v}(0,1/q^8;q^{16})$ are quasimock theta
functions, while the series $\mathcal{N}_{2v}(1,1;q)$ are
quasimodular forms. These four specializations correspond to four
important ``rank" functions (see Section 2).

\begin{theorem} \label{main2}
The functions $\mathcal{N}_{2v}(1,1;q)$ are quasimodular forms.
\end{theorem}

\begin{theorem} \label{main1}
The functions $\mathcal{N}_{2v}(1,0;q)$,
$\mathcal{N}_{2v}(1,1/q;q^2)$, and
$q^{-1}\mathcal{N}_{2v}(1,1/q^8;q^{16})$ are quasimock theta
functions.
\end{theorem}

Let us now recall what it means to be a quasimock theta function.
If $k\in \frac{1}{2}\Z\setminus \Z$, $z=x+iy$ with $x, y\in \R$,
then the weight $k$ hyperbolic Laplacian is given by
\begin{equation*}\label{laplacian}
\Delta_k := -y^2\left( \frac{\partial^2}{\partial x^2} +
\frac{\partial^2}{\partial y^2}\right) + iky\left(
\frac{\partial}{\partial x}+i \frac{\partial}{\partial y}\right).
\end{equation*}
If $v$ is odd, then define $\epsilon_v$ by
\begin{equation*}
\epsilon_v:=\begin{cases} 1 \ \ \ \ &{\text {\rm if}}\ v\equiv
1\pmod 4,\\
i \ \ \ \ &{\text {\rm if}}\ v\equiv 3\pmod 4. \end{cases}
\end{equation*}
Moreover we let $\chi$ be a Dirichlet character.
 A {\it (harmonic) weak Maass form of weight $k$ with Nebentypus $\chi$ on a subgroup
$\Gamma \subset \Gamma_0(4)$} is any smooth function $g:\H\to \C$
satisfying the following:
\begin{enumerate}
\item For all $A= \left(\begin{smallmatrix}a&b\\c&d
\end{smallmatrix} \right)\in \Gamma$ and all $z\in \H$, we
have
\begin{displaymath}
g(Az)= \leg{c}{d}\epsilon_d^{-2k} \chi(d)\,(cz+d)^{k}\ g(z).
\end{displaymath}
\item We  have that $\Delta_k g=0$. \item The function $g(z)$ has
at most linear exponential growth at all the cusps of $\Gamma$.
\end{enumerate}

In light of recent work of the first author  and Ono \cite{BO2,
Br-On1} combined with work of  Zwegers \cite{Zw1}, we now know
that what Ramanujan called mock theta functions in his last letter
to Hardy \cite{Ra} are actually ``holomorphic parts'' of weak
Maass forms.  In turn the holomorphic part of any weak Maass form
may be called a mock theta function.
In analogy with quasimodular forms \cite{kz}, a \textit{quasiweak
Maass form}  is defined to be any linear combination of
derivatives of weak Maass forms.  A function $f(q)$ is called a
\textit{quasimock theta function} if there is a quasimodular form
$h(q)$ such that $f(q) + h(q)$ is a linear combination of
derivatives of holomorphic parts of weak Maass forms.  As usual,
$q:=e^{2\pi i z}$.  Notice that taking derivatives preserves the
space of quasimock theta functions.

Our approach to Theorem \ref{main1} highlights the role that
$q$-series identities can play in the study of weak Maass forms.
Typically (see \cite{Br1}, for example) one requires lengthy and
delicate analytic calculations to determine transformation
properties.  However, we shall use $q$-series identities to
circumvent these calculations. We proceed in the same manner for
each of the three cases.  First, we use a generalized Lambert
series identity to establish the case $v=1$ by relating the
relevant function to known weak Maass forms studied in
\cite{Br-Lo2,Br-Lo3}.  Then, following the lead of
\cite{Br-Ga-Ma1}, we prove a partial differential equation and
use it to establish the case $v \geq 2$ by induction.

The functions $\mathcal{N}_{2v}(1,0;q)$ are treated in Section 4,
the $\mathcal{N}_{2v}(1,1/q;q^2)$ in Section 5, and the\\
$q^{-1}\mathcal{N}_{2v}(0,1/q^8;q^{16})$ in Section 6.  The proof
of Theorem \ref{main2} is more straightforward, following from a
certain infinite product associated with overpartition pairs. This
is discussed in Section 3.

Finally in Section 7 we take a closer look at the case $v=1$ of
\eqref{eta2vofrsn}.  Andrews observed that the function
$Spt(0,0;q)$, where
\begin{equation} \label{Sptofdeq}
Spt(d,e;q) := \frac{(-dq,-eq)_{\infty}}{(q,deq)_{\infty}}\sum_{n
\geq 1} \frac{nq^n}{(1-q^n)} - \mathcal{N}_2(d,e;q),
\end{equation}
has an elegant combinatorial interpretation and satisfies some
nice congruence properties \cite{An2}.  Further congruence
properties were found by Folsom-Ono \cite{Fo-On1} and Garvan
\cite{Ga1}.  Again these ultimately arise from the fact that
$q^{-1}Spt(0,0;q^{24})$ is a quasimock theta function.  We shall
give a combinatorial interpretation of $Spt(d,e;q)$ that reduces
to Andrews' when $d=e=0$.  Since specializations of the first term
on the right hand side of \eqref{Sptofdeq} are quasimodular forms,
the following corollary is immediate from Theorem \ref{main1}:

\begin{corollary} \label{sptquasimock}
The series $Spt(1,0;q)$, $Spt(1,1/q;q^2)$, and
$q^{-1}Spt(0,1/q^8;q^{16})$ are quasimock theta functions.
\end{corollary}
\begin{remark} \label{casede1}
As for the case $d=e=1$, it turns out that $Spt(1,1;q)$ easily
simplifies.  We have
$$
Spt(1,1;q) = -1/4 + (-q)_{\infty}^2/4(q)_{\infty}^2,
$$
which is (essentially) a modular form.
\end{remark}

Theorems \ref{main2} and \ref{main1} along with Corollary
\ref{sptquasimock} provide the theoretical framework necessary to
prove any specific number-theoretic fact about these functions.
The types of results obtainable and the methods to be employed are
well-documented in \cite{An1,An2,Br1,Br-Ga-Ma1,Fo-On1,Ga1}, and so
we shall not pursue this here.

\section{Rank moments and marked Durfee symbols}
\label{CombSection}
\subsection{A generalized $k$th symmetrized rank moment}
Recall that Dyson's \textit{rank} of a partition is the largest
part minus the number of parts \cite{dyson}.  Atkin and Garvan
\cite{atkgar} initiated the study of rank moments, the \emph{$k$th
rank moment} $N_k(n)$ being defined by
\begin{equation*} 
N_k(n) := \sum_{m \in \mathbb{Z}} m^kN(m,n).
\end{equation*}
Here $N(m,n)$ denotes the number of partitions of $n$ with rank
$m$.  Following their lead, Andrews \cite{An1} defined the
\textit{$k$th symmetrized rank moment} by
\begin{equation*} \label{Anmoment}
\eta_k(n) := \sum_{m \in \mathbb{Z}} {m + \lfloor \frac{k-1}{2}
\rfloor \choose k} N(m,n).
\end{equation*}
Evidently, the symmetrized rank moments can be expressed in terms
of the ordinary rank moments, and vice versa.  One reason to
consider the symmetrized rank moment is its simple generating
function. Namely, we have
$$
\sum_{n \geq 0}\eta_{2v}(n)q^n = \mathcal{N}_{2v}(0,0;q).
$$
When $k$ is odd, the relation $N(m,n) = N(-m,n)$ implies that both
of the moments $\eta_k(n)$ and $N_k(n)$ are $0$.

Here we will interpret the series $\mathcal{N}_{2v}(d,e;q)$ in
terms of rank moments as well, but using the rank of an
overpartition pair \cite{Br-Lo1,Lo-Ma1}. Recall that an
\textit{overpartition} $\lambda$ of $n$ is a partition of $n$ in
which the first occurrence of a number may be overlined.  An
\textit{overpartition pair} $(\lambda,\mu)$ of $n$ is a pair of
overpartitions where the sum of all of the parts is $n$. To define
the rank of an overpartition pair we use the notations
$\ell(\cdot)$ and $n(\cdot)$ for the largest part and the number
of parts of an object.  Overlining these functions indicates that
we are only considering the overlined parts.  We order the parts
of $(\lambda,\mu)$ by stipulating that for a number $k$,
\begin{equation*} \label{order}
\overline{k}_{\lambda} > k_{\lambda} > \overline{k}_{\mu} >
k_{\mu},
\end{equation*}
where the subscript indicates to which of the two overpartitions
the part belongs.
The \textit{rank of an overpartition pair} $(\lambda,\mu)$ is
\begin{equation*} \label{rankdef}
\ell((\lambda,\mu)) - n(\lambda) - \overline{n}(\mu) -
\chi((\lambda,\mu)),
\end{equation*}
where $\chi((\lambda,\mu))$ is defined to be $1$ if the largest
part of $(\lambda,\mu)$ is non-overlined and in $\mu$, and $0$
otherwise.
For example, the rank of the overpartition pair
$((\overline{6},6,5,4,4,4,\overline{3},\overline{1}),(7,7,\overline{5},2,2,2))$
is $7 - 8 - 1 - 1 = - 3$, while the rank of the overpartition pair
$((4,\overline{3},3,\overline{2},1),(4,4,4,\overline{1}))$ is $4 -
5 - 1 - 0 = -2$.

Let $N(r,s,m,n)$ denote the number of overpartition pairs of $n$
having rank $m$, such that $r$ is the number of overlined parts in
$\lambda$ plus the number of non-overlined parts in $\mu$ and $s$
is the number of parts in $\mu$.  Appealing to the case $(b,q) =
(q^{1/2},q^{1/2})$ of \cite[Thm 1.2]{Lo2}, we have the generating
function
\begin{equation} \label{pairrankgf}
N(d,e,x;q) := \sum_{r,s,n \geq 0 \atop m \in \mathbb{Z}}
N(r,s,m,n)d^re^sx^mq^n = \sum_{n \geq 0}
\frac{(-1/d,-1/e)_n(deq)^n}{(xq,q/x)_n}.
\end{equation}
This includes generating functions for several important ``ranks".
When $e = 0$ and $d=1$ we recover the generating function for
Dyson's rank for overpartitions \cite{Lo1}, and when both $d$ and
$e = 0$ we recover the generating function for Dyson's rank for
partitions. When $q=q^2$, $d=1$, and $e = 1/q$, we have the
$M_2$-rank for overpartitions \cite{Lo2}, and when $q=q^2$, $d=0$,
and $e = 1/q$, we have the $M_2$-rank for partitions without
repeated odd parts \cite{Be-Ga1,Lo-Os1}. Note that the invariance
of the right hand side of \eqref{pairrankgf} under $x
\leftrightarrow 1/x$ implies that $N(r,s,m,n) = N(r,s,-m,n)$.

We are now prepared to define the general $k$th symmetrized
rank moment for overpartition pairs.  It will be useful to also
define the ordinary \emph{$k$th rank moment for overpartition
pairs}. It is
\begin{equation}
N_k(r,s,n) := \sum_{m \in \mathbb{Z}} m^kN(r,s,m,n),
\end{equation}
and we denote its generating function by $M_k(d,e;q)$,
\begin{equation}
M_k(d,e;q) := \sum_{r,s,n \geq 0} N_k(r,s,n)d^re^sq^n.
\end{equation}
The \textit{$k$th symmetrized rank moment} is
\begin{equation}
\eta_k(r,s,n) := \sum_{m \in \mathbb{Z}} {m + \lfloor \frac{k-1}{2}
\rfloor \choose k} N(r,s,m,n).
\end{equation}

\begin{theorem} \label{ksymmomentgf}
We have
\begin{equation*}
\sum_{r,s,n \geq 0} \eta_k(r,s,n)d^re^sq^n =
\begin{cases}
0& \text{if $k$ is odd}, \\
\mathcal{N}_{2v}(d,e;q) & \text{if $k = 2v$}.
\end{cases}
\end{equation*}
\end{theorem}
\begin{proof}
Since the proof is similar to the proof of \cite[Theorem 2]{An1},
we omit most of the details.  From a limiting case of the
Watson-Whipple transformation \cite{Ga-Ra1},
\begin{equation*} \label{Watson-Whipple}
\sum_{n=0}^{\infty}
\frac{\left(aq/bc,d,e\right)_n(\frac{aq}{de})^n}{\left(q,aq/b,aq/c\right)_n} =
\frac{\left(aq/d,aq/e\right)_{\infty}}{\left(aq,aq/de\right)_{\infty}}
\sum_{n=0}^{\infty}
\frac{(a,\sqrt{a}q,-\sqrt{a}q,b,c,d,e)_n(aq)^{2n}(-1)^nq^{n(n-1)/2}}
{(q,\sqrt{a},-\sqrt{a},aq/b,aq/c,aq/d,aq/e)_n(bcde)^n},
\end{equation*}
applied to \eqref{pairrankgf} with $(a,b,c,d,e) =
(1,x,1/x,-1/d,-1/e)$, one may deduce the following alternative
form for the generating function for $N(r,s,m,n)$:
\begin{equation} \label{fromww}
\sum_{r,s,n \geq 0 \atop m \in \mathbb{Z}} N(r,s,m,n) d^re^sx^mq^n
= \frac{(-dq,-eq)_{\infty}(1-x)}{(q,deq)_{\infty}} \sum_{n \in
\mathbb{Z}}
\frac{(-de)^nq^{n(n+3)/2}(-1/d,-1/e)_n}{(-dq,-eq)_n(1-xq^n)}.
\end{equation}
Here we   note the helpful relation
\begin{equation*}
(a)_{-n} = \frac{(-1)^nq^{n+1 \choose 2}}{a^n(q/a)_n}.
\end{equation*}
Next we observe that
\begin{equation*}
\sum_{r,s,n \geq 0} \eta_{2v}(r,s,n)d^re^sq^n =
\frac{1}{(2v)!}\left(\frac{\partial^{2v}}{\partial
x^{2v}}x^{v-1}N(d,e,x;q)\right)_{x=1}.
\end{equation*}
Now computing the derivatives and simplifying as in
\cite[p.41-42]{An1}, we arrive at $\mathcal{N}_{2v}(d,e;q)$.
\end{proof}
As with ordinary partitions, the symmetrized rank moments for
overpartition pairs can be expressed in terms of the ordinary rank
moments, and vice versa. In particular we note that we have
$$N_2(r,s,n) = 2\eta_2(r,s,n).$$



\subsection{A generalized $k$-marked Durfee symbol}
The second partition-theoretic object that Andrews associated to
$\mathcal{N}_{2v}(0,0;q)$ is the $k$-marked Durfee symbol
\cite[Section 4]{An1}.  Its definition is considerably more
involved than that of the $k$th symmetrized rank moment.  We start
with $k$ copies of the natural numbers $\{1_1,2_1,3_1,\dots \}$,
$\{1_2,2_2,3_2,\dots \},\dots,\{1_k,2_k,3_k,\dots \}$.  We then
form the \textit{$k$-marked Durfee symbol} as a two-rowed array with a
subscript $S$.  Each row contains a partition using these $k$
copies of the natural numbers where parts are at most $S$. The
rows need not be of equal length.  In addition we require that:
\begin{enumerate}
\item The sequence of parts and the sequence of subscripts be
non-increasing in each row, \item Each of the subscripts smaller
than $k$ occur at least once in the top row, \item If
$M_1,N_2,\dots,V_{k-2},W_{k-1}$ are the largest parts with their
 respective subscripts in the top row, then all parts in the bottom
row with subscript $1$ lie in the interval $[1,M]$, with subscript
$2$ lie in $[M,N] \dots,$ with subscript $k-1$ lie in $[V,W]$, and
with subscript $k$ lie in $[W,S]$.
\end{enumerate}
We let $n$ be the sum of $S^2$
and all of the parts in the array and we say that the Durfee
symbol is \emph{related} to $n$.
We denote by $\mathcal{D}_k(n)$
the number of $k$-marked Durfee symbols related to $n$.  Andrews
\cite{An1} has shown that $\mathcal{D}_{v+1}(n) = \eta_{2v}(n)$.

We now define a generalized \textit{$k$-marked Durfee symbol} whose
generating function will be $\mathcal{N}_{2v}(d,e;q)$. The only
difference here is that the subscript $S$ (contributing $S^2$)
will be replaced by a triple $(S,\mu,\nu)$, $\mu$ and $\nu$ being
partitions into distinct parts between $0$ and $S-1$. For such a
partition, we say that a number $k \in [0,S-1]$ is \emph{missing}
if it does not occur.  Let $r$ denote the number of missing
numbers in $\mu$ and $s$ the number of missing numbers in $\nu$.
The number $n$ to which such a Durfee symbol is related is the sum
of $S$, all of the parts in the array, and all of the parts in
$\mu$ and $\nu$. When both $\mu$ and $\nu$ are ``full" , i.e.,
$r=s=0$, we get $S^2$, the case of the ordinary $k$-marked Durfee
symbols. For example,
\begin{equation} \label{example}
\begin{pmatrix}
4_3 & 3_2 & 3_1 & 2_1 & 1_1 \\
4_3 & 4_3 & 3_2 & 3_1 & 3_1 & 1_1
\end{pmatrix}_{4,(3,2,0),(2,1)}
\end{equation}
is a $3$-marked Durfee symbol related to $n=43$, with $r=1$ and
$s=2$.

Let $\mathcal{D}_k(r,s,n)$ be the number of generalized $k$-marked
Durfee symbols described above.  Following Andrews we define $k$
ranks associated with $k$-marked Durfee symbols.  For such a
symbol $\delta$ and for each $i$ we denote the number of entries
in the top (resp. bottom) row with subscript $i$ by
$\tau_i(\delta)$ (resp. $\beta_i(\delta)$).  Then the \textit{ith rank}
of $\delta$ is defined as
\begin{equation*} \label{ithrank}
\rho_i(\delta) :=
\begin{cases}
\tau_i(\delta) - \beta_i(\delta) - 1 & \text{for $1 \leq i <
k$},\\
\tau_i(\delta) - \beta_i(\delta)& \text{for $i = k$.}
\end{cases}
\end{equation*}
For example, the Durfee symbol in \eqref{example} has all three of
its ranks equal to $-1$.

Let $\mathcal{D}_k(r,s,m_1,m_2,\dots,m_k,n)$ denote the number of
generalized $k$-marked Durfee symbols counted by
$\mathcal{D}_k(r,s,n)$ with $i$th rank equal to $m_i$.

\begin{theorem} \label{genDsymbolgf}
For $k \geq 2$ we have
\begin{equation} \label{genDsymbolgfeq}
\begin{aligned}
&\sum_{m_1,m_2,\dots,m_k \in \mathbb{Z}} \sum_{r,s,n \geq 0}
\mathcal{D}_k(r,s,m_1,m_2,\dots,m_k,n) x_1^{m_1}x_2^{m_2}\cdots
x_k^{m_k} d^re^sq^n \\
=
&\frac{(-dq,-eq)_{\infty}}{(q,deq)_{\infty}} \sum_{n \geq 1}
\frac{(-1)^{n-1}(1+q^n)(1-q^n)^2(-1/d,-1/e)_n(de)^nq^{{n \choose
2} +kn}}{(-dq,-eq)_n\prod_{j=1}^k(1-x_jq^n)(1-q^n/x_j)}.
\end{aligned}
\end{equation}
\end{theorem}

\begin{proof}
Arguing as in \cite[Proof of Thm 10]{An1} and using the fact that
$(-1/y)_n y^n$ is the generating function for partitions into
distinct parts between $0$ and $n-1$, with the exponent of $y$
counting the number of missing numbers, we have that
\begin{equation*} \label{genDsymbolmultigf}
\begin{aligned}
&\sum_{m_1,m_2,\dots,m_k \in \mathbb{Z}} \sum_{r,s,n \geq 0}
\mathcal{D}_k(r,s,m_1,m_2,\dots,m_k,n) x_1^{m_1}x_2^{m_2}\cdots
x_k^{m_k} d^re^sq^n = \\
&\sum_{n_1 > 0 \atop n_2,n_3,\dots,n_k \geq 0}
\frac{(-1/d,-1/e)_{n_1+n_2+\cdots+n_k}(de)^{n_1+n_2+\cdots+n_k}q^{(n_1+n_2+\cdots+n_k)
+ (n_1+n_2+\cdots+n_{k-1}) + \cdots +
n_1}}{(x_1q,q/x_1)_{n_1}(x_2q^{n_1},q^{n_1}/x_2)_{n_2+1}} \\
&\times
\frac{1}{(x_3q^{n_1+n_2},q^{n_1+n_2}/x_3)_{n_3+1}
\cdots(x_kq^{n_1+\cdots+n_{k-1}},q^{n_1+\cdots+n_{k-1}}/x_k)_{n_k+1}}.
\end{aligned}
\end{equation*}
Now in the $k$-fold generalization of Watson's $q$-analogue of
Whipple's theorem \cite[eq. (2.4)]{An1}, replace $k$ by $k+1$, set
$b_j = 1/c_j = x_j$ for $1\leq j \leq k$, set $a=1$, set $b_{k+1}
= -1/d$ and $c_{k+1} = -1/e$, and finally let $N \to \infty$. Then
exactly as in \cite[Proof of Thm 3]{An1} the right hand side above
may be seen to be equal to the right hand side of
\eqref{genDsymbolgfeq}.
\end{proof}

Setting each $x_j = 1$ in Theorem \ref{genDsymbolgf},
we obtain the following:
\begin{corollary}
For $v \geq 1$ we have $\mathcal{D}_{v+1}(r,s,n) =
\eta_{2v}(r,s,n)$.
\end{corollary}

\subsection{The full rank}
We now define a statistic on generalized $k$-marked Durfee
symbols, called the \emph{full rank}.  This will not be required
in the sequel, but it plays an important role in the study of
ordinary $k$-marked Durfee symbols \cite{An1,Br1,Br-Ga-Ma1} and
will certainly do so for the generalized symbols as well. The full
rank of such a symbol $\delta$ is
$$
FR(\delta) := \rho_1(\delta) + 2\rho_2(\delta) + 3\rho_3(\delta) +
\cdots + k\rho_k(\delta).
$$
Let $NF_k(r,s,m,n)$ denote the number of generalized $k$-marked
Durfee symbols counted by $\mathcal{D}_k(r,s,n)$ whose full rank
is equal to $m$.
Evidently the two-variable
generating function for $NF_k(r,s,m,n)$ is
\begin{equation*}
\sum_{m \in \mathbb{Z} \atop r,s,n \geq 0}
NF_k(r,s,m,n)d^re^sx^mq^n = R_k\left(d,e,x,x^2,\dots,x^k;q\right),
\end{equation*}
where for $k \geq 2$, $R_k(d,e,x_1,x_2,\dots,x_k;q)$ denotes the
left-hand side of \eqref{genDsymbolgfeq}.  This function\\
$R_k(d,e,x_1,x_2,\dots,x_k;q)$ for $k \geq 2$ can in fact be
expressed in terms of $N(d,e,x;q)$.
Exactly as in Theorem 7 of \cite{An1}, we can show:
\begin{theorem}
If  $x_i \not= x_j,x_j^{-1}$ for $i \not= j$ and $x_i^2 \not=1$, then we have:
\begin{equation} \label{R1}
R_k(d,e,x_1,x_2,\dots,x_k;q) = \sum_{i=1}^k
\frac{N(d,e,x_i;q)}{\prod_{j=1 \atop j \neq
i}^k\left(x_i-x_j\right)\left(1-\frac{1}{x_ix_j}\right)}.
\end{equation}
\end{theorem}
\begin{remark}
If $x_i \in \{ x_j,x_j^{-1}\}$ or $x_i^2=1$, then a relation similar to  (\ref{R1}) can be defined via analytic continuation.
\end{remark}

\section{The case $(d,e,q) = (1,1,q)$}
Here we prove Theorem \ref{main2}. We recall from \cite{Br-Lo2}
that the two-variable generating function for
$\overline{NN}(m,n)$, the number of overpartition pairs of $n$
with rank $m$, has the very special form given by
\begin{equation} \label{veryspecial}
\sum_{n \geq 1 \atop m \in \mathbb{Z}} \overline{NN}(m,n)x^mq^n =
\frac{-4x}{(1+x)^2} +
\frac{4x(-q)_{\infty}^2}{(1+x)^2\left(xq,q/ x\right)_{\infty}}.
\end{equation}
Hence the ordinary $k$th rank moment generating function
$M_k(1,1;q)$ is simply
\begin{equation} \label{expression}
\delta_x^k \left(\frac{-4x}{(1+x)^2} +
\frac{4x(-q)_{\infty}^2}{(1+x)^2\left(xq,q/x\right)_{\infty}} \right)_{x=1},
\end{equation}
where $\delta_x := x \frac{\partial}{\partial x}$.  Now this gives
that the $M_k(1,1;q)$ are quasimodular forms (see \cite[Section
4]{atkgar} for a discussion of the quasimodularity of terms like
those in \eqref{expression}).  Since the $\mathcal{N}_{2v}(1,1;q)$
can be written as linear combinations of the $M_k(1,1;q)$, this
completes the proof of Theorem \ref{main2}.

\section{The case $(d,e,q) = (1,0,q)$}
We begin by explicitly stating and proving the case $v=1$. Let
\begin{equation}
E_2(z) := 1-24\sum_{n \geq 1} \frac{nq^n}{(1-q^n)}
\end{equation}
be the usual weight $2$ (quasimodular) Eisenstein series. Define
the integral
\begin{equation*}
  \overline{NH}(z):= \frac{1}{2\sqrt{2}\pi i} \int_{-\bar z}^{i \infty}
   \frac{\eta^2(\tau)}{\eta(2\tau)(-i (\tau+z))^{\frac32}}
  \, d \tau.
\end{equation*}
\begin{theorem} \label{d1e0v1}
The function
$$\mathcal{N}_2(1,0;q) + \frac{(-q)_{\infty}}{(q)_{\infty}}\left(\frac{1}{12} + \frac{1}{6}E_2(2z)\right)
- \overline{NH}(z)$$ is a weak Maass form of weight $\frac{3}{2}$
on $\Gamma_0(16)$.
\end{theorem}

\begin{proof} First we recall an identity involving generalized
Lambert series,
\begin{equation*} \label{Lambert}
\sum_{n \in \mathbb{Z}} \frac{(-1)^{n-1}q^{n^2}(1-x)}{\left(1-xq^n\right)} +
\sum_{n \in \mathbb{Z}} \frac{(-1)^{n}q^{n^2}(1-x)}{\left(1+xq^n\right)} =
\frac{-2\left(q^2;q^2\right)_{\infty}^2}{\left(1+1/x \right)\left(x^2q^2,q^2/x^2;q^2\right)_{\infty}}.
\end{equation*}
This is the case $y=-1/x$ of \cite[eq. (4.3), corrected]{Ch1}.
Differentiating twice with respect to $x$, setting $x=1$, and
multiplying by $(-q)_{\infty}/(q)_{\infty}$ yields
\begin{equation} \label{Lambertcor}
-4\mathcal{N}_2(1,0;q) +
4\frac{(-q)_{\infty}}{(q)_{\infty}}\sum_{n \in \mathbb{Z}}
\frac{(-1)^nq^{n^2+n}}{\left(1+q^n\right)^2} =
\frac{(-q)_{\infty}}{(q)_{\infty}} -
8\frac{(-q)_{\infty}}{(q)_{\infty}}\sum_{n \geq 1}
\frac{2nq^{2n}}{\left(1-q^{2n}\right)}.
\end{equation}
Next we recall from \cite{Br-Lo2} that the function
\begin{equation} \label{rankMaass}
\overline{\mathcal{M}}(z) :=
4\frac{(-q)_{\infty}}{(q)_{\infty}}\sum_{n \in \mathbb{Z}}
\frac{(-1)^nq^{n^2+n}}{(1+q^n)^2} - 4\overline{NH}(z)
\end{equation}
is a weak Maass form of weight $\frac{3}{2}$ on $\Gamma_0(16)$.
The theorem follows by substituting for the sum above using
\eqref{Lambertcor}.
\end{proof}
Now to prove the case $v >1$, we use a partial differential
equation which is analogous to the ``rank-crank PDE" of Atkin and
Garvan \cite{atkgar}. We define $C(x;q)$ and $C^{*}(x;q)$ by
\begin{equation*}
C(x;q) := \frac{(q)_{\infty}}{\left(xq,q/x\right)_{\infty}}, \qquad
\qquad C^{*}(x;q):=\frac{C(x;q)}{(1-x)}.
\end{equation*}
For the function $N(d,e,x;q)$ in \eqref{pairrankgf}, we define
\begin{equation}\label{rstar}
{N}^{*}(d,e,x;q):= \frac{{N}(d,e,x;q)}{(1-x)}.
\end{equation}
We use the differential operator
\begin{equation*} \label{delq}
\qquad \delta_{q}:=q \frac{\partial}{\partial q}.
\end{equation*}
Furthermore let
$$J(x;q):= \left(x,q/x\right)_{\infty}.$$ We will prove the
following partial differential equations:

\begin{theorem} \label{pde} We have
\begin{equation}\label{pdestar}
x\frac{(q)_{\infty}^2}{(-q)_{\infty}} [C^{*}(x;q)]^3 J(-x;q) =
\Bigl (2(1+x) \delta_{q} + \frac{1}{2} x + x\delta_{x} +
\frac{1}{2} (1+x) \delta_{x}^2 \Bigl )N^{*}(1,0,x;q),
\end{equation}

\noindent and
\begin{equation} \label{pdenostar}
\begin{aligned}
x\frac{(q)_{\infty}^2}{(-q)_{\infty}} & [C(x;q)]^3 J(-x;q) \\
& = \Bigl (
2(1-x)^2 (1+x) \delta_{q} + x(1+x) + 2x(1-x) \delta_{x} +
\frac{1}{2} (1+x) (1-x)^2 \delta_{x}^2  \Bigr ) N(1,0,x;q). \\
\end{aligned}
\end{equation}
\end{theorem}


\begin{proof}
Define
\begin{equation*}
S_{1}(x,\zeta;q):= \sum_{n \in \mathbb{Z}} \frac{(-1)^n \zeta^{n}
q^{n^2 + n}}{1-xq^n}.
\end{equation*}
By Lemma 4.1 in \cite{loveoz}, we have
\begin{equation} \label{iden}
\begin{aligned}
S_{1}\left(x\zeta^{-1}, \zeta^{-2};q\right) + \zeta^{2} S_{1}\left(x\zeta, \zeta^{2};q\right)
- & \zeta \frac{J\left(\zeta^2;q\right) J\left(-q;q\right)}{J\left(\zeta; q\right) J\left(-\zeta;q\right)}
S_{1}\left(x,1;q\right)
\\
& = \frac{J\left(\zeta;q\right) J\left(\zeta^{2}; q\right) J\left(-x;q\right)
(q)_{\infty}^2}{J\left(-\zeta; q\right) J\left(x\zeta;q\right) J\left(x\zeta^{-1}; q\right)
J(x;q)}. \\
\end{aligned}
\end{equation}
Equation (\ref{iden}) was one of the key results used to prove
  identities for
rank differences for overpartitions in \cite{loveoz}. Let
$g(\zeta)$ denote the right side of (\ref{iden}). Note that
$g(\zeta)$ has a double zero at $1$ and that

\begin{equation*} \label{g''}
g''(1)= \frac{2(q)_{\infty}^3}{(-q)_{\infty}^2} \left[C^{*}(x;q)\right]^3
J\left(-x;q\right).
\end{equation*}
We let $h(\zeta)$ be the sum of the first two terms on the left
side of (\ref{iden}). We  find that
$h''(1)$ equals
\begin{equation*} \label{h''}
  \sum_{n \in \mathbb{Z}} (-1)^n q^{n(n+1)} \Biggl ( \frac{8n^2 + 8n + 2}{1-xq^n} +
  \frac{2x(3+ 4n)q^n}{(1-xq^n)^2} + \frac{4x^2 q^{2n}}{(1-xq^n)^3} \Biggr )
 = \left(2+ 8\delta_{q} + 4\delta_{x} + 2\delta_{x}^2\right) \, S_{1}\left(x,1;q\right).
\end{equation*}

Letting $j(\zeta)$ be the third term on the left side of
(\ref{iden}), one can show that

\begin{equation*} \label{j''}
\begin{aligned}
j''(1) &=-4 \Biggl ( -\sum_{n=1}^{\infty} \frac{q^n}{(1+q^n)^2} -
3 \sum_{n=1}^{\infty} \frac{q^n}{(1-q^n)^2} \Biggr) S_{1}(x,1;q)
\\
&= -4 \Biggl ( -\sum_{n=1}^{\infty} \frac{q^n}{(1+q^n)^2} - 3
\Phi_{1}(q) \Biggr) S_{1}\left(x,1;q\right),
\end{aligned}
\end{equation*}

\noindent where

$$
\Phi_{1}(q):=\sum_{n=1}^{\infty} \frac{nq^n}{1-q^n}.
$$
One can check that
\begin{equation} \label{delqid}
\delta_{q} (q)_{\infty} = -\Phi_{1}(q) (q)_{\infty}.
\end{equation}

We next need the following identity, which follows from
\eqref{fromww} and the fact that
\begin{equation} \label{littlesum}
\sum_{n \in \mathbb{Z}} \frac{(-1)^nq^{n^2+n}}{(1+q^n)} =
\frac{1}{2}\frac{(q)_{\infty}}{(-q)_{\infty}}:
\end{equation}
\begin{equation} \label{rel}
xS_{1}\left(x,1;q\right)= \frac{(q)_{\infty}}{2(-q)_{\infty}} \Biggl(-1 + (1+x)
N^{*}\left(1,0,x;q\right)\Biggr).
\end{equation}

\noindent Applying $\delta_{q}$ to both sides of \eqref{rel} and
using \eqref{delqid}, we get

\begin{equation*} \label{delqs}
x\delta_{q} S_{1}\left(x,1;q\right)= \Biggl( -\Phi_{1}(q) - \sum_{k=1}^{\infty}
\frac{kq^k}{1+q^k} \Biggr)xS_{1}\left(x,1;q\right) +
\frac{(q)_{\infty}}{2(-q)_{\infty}} \delta_{q} \left(1+x\right)
N^{*}\left(1,0,x;q\right).
\end{equation*}

\noindent Similarly, we find that

\begin{equation*} \label{delzs}
x\delta_{x} S_{1}\left(x,1;q\right) = \frac{(q)_{\infty}}{2(-q)_{\infty}}
\delta_{x} \left(1+x\right) N^{*}\left(1,0,x;q\right) - xS_{1}\left(x,1;q\right)
\end{equation*}
and
\begin{equation*}
\label{delz2s} x\delta_{x}^2 S_{1}\left(x,1;q\right) = xS_{1}\left(x,1;q\right) +
\frac{(q)_{\infty}}{2(-q)_{\infty}} \left(\delta_{x}^2 -
2\delta_{x}\right) N^{*}\left(1,0,x;q\right) \left(1+x\right).
\end{equation*}
Combining the above now easily yields

\begin{multline}
 \label{ans}
x\frac{(q)_{\infty}^2}{(-q)_{\infty}} \left[C^{*}(x;q)\right]^3 J(-x;q) 
=
\Bigl (2\delta_{q} + \frac{1}{2} \delta_{x}^2 \Bigr )N^{*}\left(1,0,x;q\right) \left(1+x\right) \\
+ 2x\frac{(-q)_{\infty}}{(q)_{\infty}} S_{1}\left(x,1;q\right) \Biggl [
\sum_{n=1}^{\infty} \frac{nq^n}{1-q^n} - 2\sum_{n=1}^{\infty}
\frac{nq^n}{1+q^n} +  \sum_{n=1}^{\infty} \frac{q^n}{\left(1+q^n\right)^2}
\Biggr ].
\end{multline}

\noindent Note that the terms in brackets in (\ref{ans}) sum to
$0$.  This may be seen, for example, by writing these sums in
terms of divisor functions.

Next an application of the product rule yields (\ref{pdestar}).
From (\ref{rstar}) we find that

\begin{equation} \label{e1}
\delta_{x} N^{*}\left(1,0,x;q\right)=\frac{\delta_{x} N\left(1,0,x;q\right) +
xN^{*}\left(1,0,x;q\right)}{1-x},
\end{equation}

\begin{equation} \label{e2}
\delta_{x}^{2} N^{*}\left(1,0,x;q\right)=\frac{\delta_{x}^{2} N\left(1,0,x;q\right) +
2x\delta_{x} N^{*}\left(1,0,x;q\right) + xN^{*}\left(1,0,x;q\right)}{1-x},
\end{equation}

\begin{equation} \label{e3}
\delta_{q} N^{*}\left(1,0,x;q\right)=\frac{\delta_{q} N\left(1,0,x;q\right)}{1-x}.
\end{equation}
This easily yields
%
(\ref{pdenostar}).
\end{proof}

We may now prove the case $v > 1$ inductively using Theorem
\ref{pde}. Actually we shall argue using the ordinary rank moment
generating functions $M_{2v}(1,0;q)$, but as we have already
mentioned the functions $\mathcal{N}_{2v}(1,0;q)$ may be written
in terms of the $M_{2v}(1,0;q)$.  First, since $M_2(1,0;q) =
2\mathcal{N}_2(1,0;q)$, Theorem \ref{d1e0v1} implies that the
former is a quasimock theta function. Next, appealing to
\eqref{pairrankgf}, we have that the generating function for
$M_v(1,0;q)$ is $\left[\delta_x^v N(1,0,x;q) \right]_{x=1}$.
Now apply $\delta_x$ to (\ref{pdenostar}) $2v$ times and then set
$x=1$. We first consider the left-hand side. First, from
\cite[Section 4]{atkgar} we have that $\left[\delta_x^r C(x;q)
\right]_{x=1}$ is a quasimodular form. Moreover we have
\begin{eqnarray*}
\delta_x \left( J(-x;q) \right)
&=& \left( \frac{x}{1+x} +
x \sum_{m=1}^{\infty} \frac{q^m}{1+xq^m} -
x^{-1} \sum_{m=1}^{\infty} \frac{q^m}{1+x^{-1}q^m} \right) J(-x;q)\\
&=&
\left(  \frac{x}{x+1}-\sum_{m=1}^{\infty}  \frac{(-1)^mq^m}{1-q^m} \left(x^m-x^{-m} \right)  \right)
J(-x;q).
\end{eqnarray*}
Thus $\left[\delta_x^{2v} J(-x;q)\right]_{x=1}$ is a linear
combination of terms of the form
\begin{equation*}
 \left( c_r - \left(1-(-1)^{r} \right)
 \sum_{m=1}^{\infty}  \frac{(-1)^m\, m^r\, q^m}{1-q^m} \right)^l
J(-1;q)
\end{equation*}
for integers $r$ and $l$ with some constant $c_r$. The theory of
Eisenstein series on congruence subgroups (see Section III.3 in
\cite{Kob}) yields that up to a constant term the above sum is a
modular form for odd $r \geq 3$.  For $r=1$, observe that
\begin{equation*}
\delta_q \frac{(q^2;q^2)_{\infty}}{(q;q^2)_{\infty}} =
\frac{-\left(q^2;q^2\right)_{\infty}}{\left(q;q^2\right)_{\infty}} \sum_{n \geq 1}
\frac{(-1)^nn\, q^n}{1-q^n},
\end{equation*}
and hence we have a quasimodular form in this case.  So, applying
$\delta_x$ $2v$ times and then setting $x=1$ yields a quasimodular
form on the left hand side.  Now the claim follows by induction,
since the  term $\delta_x^{2v}N(1,0,x;q)$  occurs with
multiplicity $(2v-2)(2v-1) \not= 0$ for $v>1$ and the other terms
on the right are derivatives of quasimock theta functions.

\section{The case $(d,e,q) = (1,1/q,q^2)$}
We begin again with the case $v=1$.
\begin{theorem} \label{d1e1qv1}
The function
\begin{equation*}
\mathcal{N}_2\left(1,1/q;q^2\right) +
\frac{(-q)_{\infty}}{(q)_{\infty}}\left(\frac{1}{12} +
\frac{1}{24}E_2\left(2z\right)\right) - \frac{\overline{NH}(z)}{2}
\end{equation*}
is a weak Maass form of weight $3/2$ on $\Gamma_0(16)$.
\end{theorem}

\begin{proof}
We begin with the following identity, which is obvious:
\begin{equation*}
\sum_{n \geq 1} \frac{(-1)^nq^{n^2+n}\left(1+q^{2n}\right)}{\left(1-q^{2n}\right)^2} +
2\sum_{n \geq 1} \frac{(-1)^nq^{n^2+2n}}{\left(1-q^{2n}\right)^2} = \sum_{n
\geq 1} \frac{(-1)^nq^{n^2+n}}{\left(1-q^n\right)^2}.
\end{equation*}
In terms of symmetrized rank moments, this says that
\begin{equation} \label{into}
\frac{(-q)_{\infty}}{(q)_{\infty}}\sum_{n \geq 1}
\frac{(-1)^nq^{n^2+n}\left(1+q^{2n}\right)}{\left(1-q^{2n}\right)^2} -
\mathcal{N}_2\left(1,1/q;q^2\right) = -\frac{\mathcal{N}_2\left(1,0;q\right)}{2}.
\end{equation}
Now, taking $(a,b,c,d,q) = (1,b,1/b,\infty,q^2)$ in the $_6\phi_5$
summation \cite{Ga-Ra1},
\begin{equation*}
\sum_{n \geq 0}
\frac{(1-aq^{2n})(a,b,c,d)_n\left(aq/bcd\right)^n}{(1-a)\left(q,aq/b,aq/c,aq/d\right)_n}
=
\frac{\left(aq,aq/bc,aq/bd,aq/cd\right)_{\infty}}{\left(aq/b,aq/c,aq/d,aq/bcd\right)_{\infty}},
\end{equation*}
we obtain
\begin{equation*}
1 + \sum_{n \geq 1}
\frac{(1+q^{2n})\left(b,1/b;q^2\right)_n(-1)^nq^{n^2+n}}{\left(bq^2,q^2/b;q^2\right)_n}
= \frac{\left(q^2;q^2\right)_{\infty}^2}{\left(bq^2,q^2/b;q^2\right)_{\infty}}.
\end{equation*}
Taking $\partial^2/\partial b^2$, setting $b=1$, and substituting
into \eqref{into} gives
\begin{equation*}
\frac{(-q)_{\infty}}{(q)_{\infty}}\sum_{n \geq
1}\frac{nq^{2n}}{\left(1-q^{2n}\right)} + \mathcal{N}_2\left(1,1/q;q^2\right) =
\frac{\mathcal{N}_2\left(1,0;q\right)}{2}.
\end{equation*}
Applying Theorem \ref{d1e0v1} completes the proof.
\end{proof}

Next we deduce a PDE for $N(1,1/q,x;q^2)$.

\begin{theorem} \label{pde2}
We have
\begin{multline} \label{pdeoverstar}
2x \left(q^2; q^2\right)_{\infty}^2  \left[C^{*}(x;q^2)\right]^3 J\left(-x; q\right)= \Bigl (
(1+x)\delta_{q} + x + 2x\delta_{x} + (1+x)\delta_{x}^2 \Bigr)
N^{*}\left(1,1/q,x; q^2\right),
\end{multline}

\noindent and

\begin{equation} \label{pdeovernostar}
\begin{aligned}
2x & \left(q^2; q^2\right)_{\infty}^2 \left[C(x;q^2)\right]^3 J(-x; q) \\
& = \Bigl
((1+x)(1-x)^2 \delta_{q} + 2x(1+x) + 4x(1-x)\delta_{x} +
(1+x)(1-x)^2 \delta_{x}^2 \Bigr) N\left(1,1/q,x; q^2\right). \\
\end{aligned}
\end{equation}
\end{theorem}

\begin{proof}
The proof is similar to that of Theorem
\ref{pde}. We define

\begin{equation*}
S_{2}\left(x,\zeta;q\right):= \sum_{n \in \mathbb{Z}} \frac{(-1)^n \zeta^{n}
q^{n^2 + 2n}}{1-xq^{2n}}.
\end{equation*}

Taking $(a_1,a_2,b_1,b_2,b_3,q) = (-x,-xq,x\zeta,x/\zeta,x,q^2)$
in the case $(r,s) = (2,3)$ of Theorem 2.1 of \cite{Ch1}, we
obtain

\begin{equation} \label{iden2}
\begin{aligned}
S_{2}\left(x\zeta^{-1}, \zeta^{-1};q\right) + \zeta^{2} S_{2}\left(x\zeta, \zeta; q\right) +
&  2
\frac{J\left(\zeta^2;q^2\right) (-q; q)_{\infty}^2}{J\left(-\zeta; q\right)
J\left(\zeta^{-1};q^2\right)} S_{2}\left(x,1;q\right) \\
&
 = \frac{J(-x;q) J\left(\zeta^{2}; q^2\right)
J\left(\zeta; q^2\right) \left(q^2; q^2\right)_{\infty}^2}{J\left(x\zeta; q^2\right)
J\left(x\zeta^{-1};q^2\right) J\left(-\zeta; q\right) J\left(x;q^2\right)}. \\
\end{aligned}
\end{equation}

Let $g(\zeta)$ denote the right side of (\ref{iden2}). Note that
$g(\zeta)$ has a double zero at $1$ and that

\begin{equation*} \label{overg''}
g''(1)= \frac{2\left(q^2; q^2\right)_{\infty}^3}{(-q; q)_{\infty}^2}
\left[C^{*}(x;q^2)\right]^3 J(-x;q).
\end{equation*}
We let $h(\zeta)$ be the sum of the first two terms on the left
side of (\ref{iden2}). We find that $h''(1)$ equals

\begin{equation*} \label{overh''}
  \sum_{n \in \mathbb{Z}} (-1)^n q^{n(n+2)} \Biggl ( \frac{2n^2 + 4n + 2}{1-xq^{2n}} +
  \frac{2x(3+ 2n)q^{2n}}{(1-xq^{2n})^2} + \frac{4x^2 q^{4n}}{(1-xq^{2n})^3} \Biggr )
 = \left(2+ 2\delta_{q} + 4\delta_{x} + 2\delta_{x}^2 \right) \, S_{2}\left(x,1;q\right).
\end{equation*}

Letting $j(\zeta)$ be the third term on the left side of
(\ref{iden2}), one can show that

\begin{equation*} \label{overj''}
j''(1) =4\Biggl (\sum_{n=1}^{\infty} \frac{q^{n}}{(1+q^{n})^2} + 3
\Phi_{1}\left(q^2\right) \Biggr) S_{2}\left(x,1;q\right).
\end{equation*}

 We next need the following identity, which is again a consequence
 of \eqref{fromww} and \eqref{littlesum}:
\begin{equation} \label{overrel}
xS_{2}\left(x,1;q\right)= \frac{(q)_{\infty}}{2(-q)_{\infty}} \left(-1 +
(1+x)N^{*}\left(1,1/q,x;q^2\right)\right).
\end{equation}

\noindent Applying $\delta_{q}$ to both sides of (\ref{overrel})
and using (\ref{delqid}), we get

\begin{equation*} \label{deloverqs}
x\delta_{q} S_{2}\left(x,1;q\right)= \Biggl( -\Phi_{1}(q) - \sum_{k=1}^{\infty}
\frac{kq^{k}}{1+q^{k}} \Biggr)xS_{2}\left(x,1;q\right) +
\frac{(q)_{\infty}}{2(-q)_{\infty}} \delta_{q}
(1+x)N^{*}\left(1,1/q,x;q^2\right).
\end{equation*}

\noindent Similarly, we find that

\begin{eqnarray*} \label{deloverzs}
x\delta_{x} S_{2}\left(x,1;q\right) = \frac{(q)_{\infty}}{2(-q)_{\infty}}
\delta_{x} (1+x)N^{*}\left(1,1/q,x;q^2\right) - xS_{2}\left(x,1;q\right),\\
\label{deloverz2s} x\delta_{x}^2 S_{2}\left(x,1;q\right) = xS_{2}\left(x,1;q\right) +
\frac{(q)_{\infty}}{2(-q)_{\infty}} \left(\delta_{x}^2 -
2\delta_{x}\right) N^{*}\left(1,1/q,x;q^2\right)(1+x).
\end{eqnarray*}
Combining the above and simplifying now yields

\begin{multline}
 \label{overans}
2x \left(q^2; q^2\right)_{\infty}^2 \left[C^{*}(x;q^2)\right]^3 
J(-x;q) = \Bigl
((1+x)\delta_{q} + x + 2x\delta_{x} + (1+x)\delta_{x}^2 \Bigr )
N^{*}\left(1,1/q,x;q^2\right) \\
 + 2x \frac{(-q)_{\infty}}{(q)_{\infty}} S_{2}\left(x,1;q\right) \Biggl [
-\Phi_{1}(q) - \sum_{n=1}^{\infty} \frac{nq^n}{1+q^n} + 2 \sum_{n=1}^{\infty} \frac{q^n}{\left(1+ q^{n}\right)^2} + 6\Phi_{1}\left(q^2\right) \Biggr ]. \\
\end{multline}

\noindent Note that the terms in brackets in (\ref{overans}) sum to 0 and thus we have (\ref{pdeoverstar}). Using the analogues of
equations (\ref{e1}), (\ref{e2}), and (\ref{e3}), we may obtain
(\ref{pdeovernostar}).
\end{proof}

Now the general case $v >1$ follows just as for
$\mathcal{N}_{2v}(1,0;q)$ in the previous section.  We omit the
details.

\section{The case $(d,e,q) = (0,1/q,q^2)$}
Let us begin again with the case $v=1$.  Before stating it, we
need a lemma. Define the function $g(z)$ by
\begin{equation*}
g(z):= \frac{\left( q;q^2\right)_{\infty}}{\left(
q^2;q^2\right)_{\infty}} \sum_{n \in \Z} \frac{q^{
2n^2+3n+1}}{\left(1-q^{2n+1} \right)^2}
\end{equation*}
and the integral $NH_2(z)$ by
\begin{equation*}
NH_2(z) :=  \frac{1}{4 \sqrt{2} i \pi} \int_{- \bar{z}}^{i \infty}
\frac{\eta^2(16 \tau)}{\eta(8 \tau) \left(- i \left( \tau+z
\right) \right)^{\frac32 }} d \tau.
\end{equation*}
Moreover let
\begin{equation*}
\overline{\mathcal{M}}_2(z):= q^{-1}g\left(8z\right) - NH_2(z).
\end{equation*}
\begin{lemma} \label{another}
The function $\overline{\mathcal{M}}_2(z)$ is a weak Maass form of
weight $\frac32$ on $\Gamma_0(16)$.
\end{lemma}
\begin{proof}
This will follow from the work in \cite[Section 4]{Br-Lo2}. First,
recall that the function $\overline{\mathcal{M}}(z)$ defined in
 $(\ref{rankMaass})$ is a weak Maass form of weight $\frac32$
 on $\Gamma_0(16)$.
 As in the case of classical modular forms one can show
 that the function
 \begin{equation*}
 \overline{\mathcal{N}}(z):=
 \frac{1}{2 \sqrt{2}} \left( - i 16z \right)^{-\frac32} \overline{\mathcal{M}} \left(-\frac{1}{16z} \right)
 \end{equation*}
 is also a weak Maass form of weight $\frac32$ on $\Gamma_0(16)$.
 It turns out that $\overline{\mathcal{N}}(z) = \overline{\mathcal{M}}_2(z)$.
 To see this, observe that the transformation law for $\overline{\mathcal{M}}(z)$
 (see Corollary 4.4 and Lemma 4.5 of \cite{Br-Lo2}) implies that
 \begin{equation*}
 \mathcal{M} \left( -\frac{1}{z}\right)
 = 2 \sqrt{2} \left(- i z \right)^{ \frac32} \mathcal{U} \left( \frac{z}{2}\right)
 - \frac{2}{ \pi i} (- i z)^{\frac32 } \int_{- \bar{z}}^{i \infty}
 \frac{\eta^2(\tau)}{\eta \left(\frac{\tau}{2}\right)
  \left(- i  \left( \tau +z \right)\right)^{\frac32}
} d \tau,
 \end{equation*}
 where
 \begin{equation*}
 \mathcal{U}(z):= \frac{ \eta(z) }{ \eta^2(2z)  }
 \sum_{  \substack{  n \in \Z  \\ n \text{ odd}     }  }
 \frac{q^{ \frac{n}{2} (n+1)}}{\left( 1-q^n\right)^2}.
 \end{equation*}
 To finish the proof we make
  the change of variables $z \mapsto 16 z$ and observe that $\mathcal{U}(z)= q^{-\frac18}g(z)$.
  \end{proof}
 \begin{remark}
 We note that there is a typo in the definition of the function
 $\mathcal{U}(z)$ just above Corollary 4.2 in \cite{Br-Lo2}.
 \end{remark}

\begin{theorem} \label{d0e1qv1}
The function
\begin{equation*}
q^{-1}\mathcal{N}_2\left(0,-1/q^8;q^{16}\right) -
\frac{\eta(8z)}{24\eta^2(16z)}\left(1-E_2\left(8z\right)\right) + NH_2(z)
\end{equation*}
is a weak Maass form of weight $3/2$ on $\Gamma_0(16)$.
\end{theorem}
\begin{proof}
Here we use the identity
\begin{equation*}
\sum_{n \in \mathbb{Z}}
\frac{(-1)^{n+1}q^{2n^2-n}}{\left(1-xq^{2n}\right)} + \sum_{n \in
\mathbb{Z}} \frac{(-1)^{n}q^{2n^2+n}}{\left(1+xq^{2n+1}\right)} =
\frac{\left(-q,q^2;q^2\right)_{\infty}^2}{\left(1/x,xq^2,-xq,-q/x;q^2\right)_{\infty}},
\end{equation*}
which is the case $q=q^2$ and $y = -1/xq$ of \cite[eq. (4.3),
corrected]{Ch1}.  Differentiating twice with respect to $x$,
setting $x=1$, and multiplying by
$(q;q^2)_{\infty}/(q^2;q^2)_{\infty}$ gives
\begin{equation*}
\frac{-\left(q;q^2\right)_{\infty}}{\left(q^2;q^2\right)_{\infty}}\sum_{n \in \mathbb{Z}
\setminus \{ 0 \}} \frac{q^{2n^2+n}}{\left(1-q^{2n}\right)^2} +
\frac{\left(q;q^2\right)_{\infty}}{\left(q^2;q^2\right)_{\infty}}\sum_{n \in \mathbb{Z}}
\frac{q^{2n^2+3n+1}}{\left(1-q^{2n+1}\right)^2} =
\frac{\left(q;q^2\right)_{\infty}}{\left(q^2;q^2\right)_{\infty}} \sum_{n \geq 1}
\frac{nq^n}{\left(1-q^n\right)}.
\end{equation*}
Notice that the first term on the left hand side above is
$\mathcal{N}_2(0,-1/q;q^2)$ and the second is $g(z)$.  The theorem
now follows from Lemma \ref{another}.
\end{proof}
\begin{remark}
Notice that by replacing $z$ by $z+\frac{1}{16}$ (i.e. $q$ by
$e^{\pi i/8}q$), we have that
$q^{-1}\mathcal{N}_2(0,-1/q^8;q^{16})$ is also a quasimock theta
function, although the corresponding weak Maass form is on a much
smaller group ( $\Gamma_1(256)$, for example).
\end{remark}

 We now prove a PDE for $N(0,1/q,x;q^2)$.
\begin{theorem} \label{pde3}
We have
\begin{equation} \label{pdem2star}
2x \frac{(q^2; q^2)_{\infty}^2}{(-q; q^2)_{\infty}}
\left[C^{*}(x;q^2)\right]^3 J\left(-xq; q^2\right) = \Bigl ( 2\delta_{q} + \delta_{x} +
\delta_{x}^2 \Bigr ) N^{*}\left(0,1/q,x; q^2\right),
\end{equation}

\noindent and

\begin{equation} \label{pdem2nostar}
\begin{aligned}
2x \frac{\left(q^2; q^2\right)_{\infty}^2}{\left(-q; q^2\right)_{\infty}} & \left[C\left(x; q^2\right)\right]^3
J\left(-xq; q^2\right) \\
& = \Bigl (2(1-x)^2 \delta_{q} + (1+x)(1-x) \delta_{x} +
2x + (1-x)^2 \delta_{x}^2 \Bigr) N\left(0,1/q,x; q^2\right). \\
\end{aligned}
\end{equation}
\end{theorem}

\begin{proof}
To prove (\ref{pdem2star}), we first define

\begin{equation*}
S_{3}\left(x,\zeta;q\right):= \sum_{n \in \mathbb{Z}} \frac{(-1)^n \zeta^{n}
q^{2n^2 + 3n}}{1-xq^{2n}}.
\end{equation*}

By Lemma 4.1 in \cite{Lo-Os1}, we have

\begin{equation} \label{iden1}
\begin{aligned}
S_{3}\left(x\zeta^{-1}, \zeta^{-2};q\right) + \zeta^{3} S_{3}\left(x\zeta, \zeta^{2};q\right) -
& \zeta \frac{J\left(\zeta^2;q^2\right) \left(-q; q^2\right)_{\infty}^2}{J\left(\zeta; q^2\right)
J\left(-q\zeta;q^2\right)} S_{3}\left(x,1;q\right) \\
& = \frac{J\left(-xq;q^2\right) J\left(\zeta^{2}; q^2\right)
J\left(\zeta; q^2\right) \left(q^2; q^2\right)_{\infty}^2}{J\left(x\zeta^{-1}; q^2\right)
J\left(x\zeta;q^2\right) J\left(-q\zeta; q^2\right) J\left(x;q^2\right)}. \\
\end{aligned}
\end{equation}

Equation (\ref{iden1}) was one of the key results used to prove
identities for
$M_2$-rank differences for partitions without repeated odd parts
in \cite{Lo-Os1}. Let $g(\zeta)$ denote the right side of
(\ref{iden1}). Note that $g(\zeta)$ has a double zero at $1$ and
that

\begin{equation*} \label{m2g''}
g''(1)= \frac{4\left(q^2; q^2\right)_{\infty}^3}{\left(-q; q^2\right)_{\infty}^2}
\left[C^{*}\left(x;q^2\right)\right]^3 J\left(-xq;q^2\right).
\end{equation*}
We let $h(\zeta)$ be the sum of the first two terms on the left
side of (\ref{iden1}). We find that $h''(1)$ equals

\begin{equation*} \label{m2h''}
  \sum_{n \in \mathbb{Z}} (-1)^n q^{n(2n+3)} \Biggl ( \frac{8n^2 + 12n + 6}{1-xq^{2n}}
  + \frac{8x(1+ n)q^{2n}}{(1-xq^{2n})^2} + \frac{4x^2 q^{4n}}{(1-xq^{2n})^3} \Biggr )
 = (6+ 4\delta_{q} + 6\delta_{x} + 2\delta_{x}^2 ) \, S_{3}\left(x,1;q\right).
\end{equation*}

Letting $j(\zeta)$ be the third term on the left side of
(\ref{iden1}), one can show that

\begin{equation*} \label{m2j''}
j''(1) =-2 \Biggl (1 - 2\sum_{n=1}^{\infty}
\frac{q^{2n-1}}{(1+q^{2n-1})^2} - 6 \Phi_{1}(q^2) \Biggr)
S_{3}\left(x,1;q\right).
\end{equation*}

One can check that
\begin{equation} \label{delq2id}
\delta_{q} \left(q^2; q^2\right)_{\infty} = -2\Phi_{1}\left(q^2\right) \left(q^2;
q^2\right)_{\infty}
\end{equation}

\noindent and

\begin{equation} \label{delqq2id}
\delta_{q} \left(-q; q^2\right)_{\infty} = \left(-q; q^2\right)_{\infty}
\sum_{k=0}^{\infty} \frac{(2k+1)q^{2k+1}}{\left(1+q^{2k+1}\right)}.
\end{equation}

From \eqref{fromww} we have the following identity:
\begin{equation} \label{m2rel}
xS_{3}\left(x,1;q\right)= \frac{\left(q^2; q^2\right)_{\infty}}{\left(-q; q^2\right)_{\infty}} \Bigl(-1
+ N^{*}\left(0,1/q,x;q^2\right)\Bigr).
\end{equation}

\noindent Applying $\delta_{q}$ to both sides of (\ref{m2rel}) and
using (\ref{delq2id}) and (\ref{delqq2id}), we get

\begin{equation*} \label{delm2qs}
x\delta_{q} S_{3}\left(x,1;q\right)= \Biggl( -2\Phi_{1}\left(q^2\right) -
\sum_{k=0}^{\infty} \frac{(2k+1)q^{2k+1}}{1+q^{2k+1}}
\Biggr)xS_{3}\left(x,1;q\right) + \frac{\left(q^2; q^2\right)_{\infty}}{\left(-q; q^2\right)_{\infty}}
\delta_{q} N^{*}\left(0,1/q,x;q^2\right).
\end{equation*}

\noindent Similarly, we find that

\begin{eqnarray*} \label{delm2zs}
x\delta_{x} S_{3}(x,1;q) &=& \frac{\left(q^2; q^2\right)_{\infty}}{(-q;
q^2)_{\infty}}
\delta_{x} N^{*}\left(0,1/q,x;q^2\right) - xS_{3}\left(x,1;q\right),\\
\label{delm2z2s} x\delta_{x}^2 S_{3}\left(x,1;q\right) &=& xS_{3}\left(x,1;q\right) +
\frac{\left(q^2;q^2\right)_{\infty}}{\left(-q; q^2\right)_{\infty}} \left(\delta_{x}^2 -
2\delta_{x}\right) N^{*}\left(0,1/q,x;q^2\right).
\end{eqnarray*}
Combining the above now yields

 \begin{equation} \label{m2ans}
\begin{aligned}
2x\frac{\left(q^2; q^2\right)_{\infty}^2}{\left(-q; q^2\right)_{\infty}}
& \left[C^{*}\left(x;q^2\right)\right]^3 J\left(-xq;q^2\right) = \Bigl (2\delta_{q} + \delta_{x} +
\delta_{x}^2 \Bigr ) N^{*}\left(0,1/q,x;q^2\right) \\
& + 4x \frac{\left(-q; q^2\right)_{\infty}}{\left(q^2; q^2\right)_{\infty}} S_{3}\left(x,1;q\right) \Biggl [ \Phi_{1}\left(q^2\right) - \sum_{n=0}^{\infty} \frac{(2n+1) q^{2n+1}}{1+q^{2n+1}} + \sum_{n=1}^{\infty} \frac{q^{2n-1}}{\left(1+ q^{2n-1}\right)^2} \Biggr ].
\end{aligned}
\end{equation}

\noindent Observe that the terms in brackets in (\ref{m2ans}) sum to 0 and thus (\ref{pdem2star}) follows.  Using the analogues of
equations (\ref{e1}), (\ref{e2}) and (\ref{e3}), we obtain
(\ref{pdem2nostar}).
\end{proof}

The case $v > 1$ follows by induction as before.

\section{The smallest parts function}

We now define a \textit{smallest parts function} $spt(r,s,n)$ in the
context of overpartition pairs $(\lambda,\mu)$ of $n$.  It is the
total number of appearances of the smallest parts in all of the
overpartition pairs $(\lambda,\mu)$ of $n$, where $r$ is the
number of overlined parts in $\lambda$ plus the number of
non-overlined parts in $\mu$ and $s$ is the number of parts in
$\mu$, such that the smallest part in $(\lambda,\mu)$ only occurs
non-overlined and only in $\lambda$. Thus overpartition pairs like
$((\overline{6},6,5,4,4,4,\overline{3},\overline{1}),(7,7,\overline{5},2,2,2))$
contribute nothing to $spt(r,s,n)$.  We have the following
generating function:

\begin{theorem} \label{sptgf}
Recalling the definition of $Spt(d,e;q)$ from the introduction, we
have
\begin{equation*}
\sum_{r,s,n \geq 0} spt(r,s,n) d^re^sq^n = Spt(d,e;q).
\end{equation*}
\end{theorem}

\begin{proof}
We proceed as in \cite[Proof of Theorem 4]{An2}.  Briefly, we have
that
\begin{equation*}
\begin{aligned}
& \sum_{r,s,n \geq 0} spt(r,s,n)d^re^sq^n =
\frac{\left(-dq,-eq\right)_{\infty}}{\left(deq,q\right)_{\infty}}\sum_{n
\geq 1}
\frac{\left(q,deq\right)_nq^n}{\left(1-q^n\right)^2\left(-dq,-eq\right)_n} \\
&= \frac{-\left(-dq,-eq\right)_{\infty}}{2\left(deq,q\right)_{\infty}}\left[
\frac{\partial ^2}{\partial x^2}\sum_{n \geq 0}
\frac{\left(deq\right)_n\left(x,1/x\right)_nq^n}{\left(q,-dq,-eq\right)_n} \right]_{x=1} \\
&=
\frac{-\left(-dq,-eq\right)_{\infty}}{2\left(deq,q\right)_{\infty}}\left[\frac{\partial
^2}{\partial x^2} \frac{\left(xq,q/x\right)_{\infty}}{(q)_{\infty}^2}\left(1+
\sum_{n \geq 1}
\frac{(-de)^nq^{n(n+3)/2}(1+q^n)\left(-1/d,-1/e,x,1/x\right)_n}{\left(q/x,xq,-dq-eq\right)_n}\right)\right]_{x=1}
\\
&= Spt(d,e;q),
\end{aligned}
\end{equation*}
where the penultimate line follows from the Watson-Whipple
transformation.
\end{proof}
As noted in the introduction, Corollary \ref{sptquasimock} follows
immediately from Theorem \ref{main1}, or more explicitly, from
Theorems \ref{d1e0v1}, \ref{d1e1qv1}, and \ref{d0e1qv1}. Regarding
Remark \ref{casede1}, it is not hard to see that when $d=e=1$, the
symmetry in the overpartition pairs means that the smallest parts
function for $n
>0$ is just counting $\overline{pp}(n)/4$, where
$\overline{pp}(n)$ is the number of overpartitions of $n$. From
the case $x=1$ of \eqref{veryspecial}, the generating function for
$\overline{pp}(n)$ is $(-q)_{\infty}^2/(q)_{\infty}^2$.









\end{document}